\documentclass[12pt]{article}

\usepackage{amsmath, epsfig, cite}
\usepackage{amssymb}
\usepackage{amsfonts}
\usepackage{latexsym}
\usepackage{epsfig,graphics,pstricks,ifpdf}
\usepackage{amsthm}

\makeatletter
\renewcommand{\@seccntformat}[1]{{\csname the#1\endcsname}.\hspace{.5em}}
\makeatother

\newtheorem{thm}{Theorem}[section]

\newtheorem{cor}[thm]{Corollary}
\newtheorem{lem}[thm]{Lemma}

\newcommand{\pf}{\noindent{\it Proof.} }

\def\Z{\mathbb{Z}}

\def\0{{\bf 0}}

\numberwithin{equation}{section}

\renewcommand{\qed}{\hfill$\Box$\medskip}

\begin{document}

\begin{center}
{\Large\bf On arithmetic partitions of $\Z_n$}
\end{center}
\vskip 2mm \centerline{Victor J. W. Guo$^1$  and Jiang Zeng$^{2}$}

\begin{center}
{\footnotesize $^1$Department of Mathematics, East China Normal
University,\\ Shanghai 200062,
 People's Republic of China\\
{\tt jwguo@math.ecnu.edu.cn,\quad http://math.ecnu.edu.cn/\textasciitilde{jwguo}}\\[10pt]
$^2$Universit\'e de Lyon; Universit\'e Lyon 1; Institut Camille
Jordan, UMR 5208 du CNRS;\\ 43, boulevard du 11 novembre 1918,
F-69622 Villeurbanne Cedex, France\\
{\tt zeng@math.univ-lyon1.fr,\quad
http://math.univ-lyon1.fr/\textasciitilde{zeng}} }
\end{center}

\vskip 0.7cm \noindent{\small{\bf Abstract.}
Generalizing
a classical problem in enumerative combinatorics,  Mansour and Sun
counted the number of subsets of $\Z_n$ without certain separations.
Chen, Wang, and Zhang then studied the  problem of
 partitioning $\Z_n$ into arithmetical progressions of
 a given type under some  technical conditions. In this paper, we
improve  on their main theorems
 by applying a convolution formula for  cyclic multinomial coefficients due to
Raney-Mohanty.}




\section{Introduction}

In his solution of \emph{probl\`eme des m\'enages} Kaplansky \cite{Kaplansky}
showed that the number of ways of selecting  $k$ elements, no two consecutive,
from $n$ objects arrayed on a cycle is $\frac{n}{n-k}{n-k\choose k}$.
Let $\Z_n:=\{0,1,\ldots, n-1\}$ be the set of congruence classes modulo $n$ with usual arithmetic.
Then Yamamoto \cite{Yamamoto} (see also \cite[p. 222]{Riordan})
proved that if $n\geqslant pk+1$ the number of ways of selecting $k$ elements from $\Z_n$,
no two consecutive,  is
\begin{align}\label{eq:yama}
\frac{n}{n-pk}{n-pk\choose k},
\end{align}
when $i\pm1,\ldots, i\pm p$ are regarded as consecutive to $i$.

In the last three decades a lot of  generalizations and variations of Kaplansky's problem
have been studied by several authors (see, for example, \cite{Konvalina,Chu,Hwang,Hwang2,HKW,KP,Moser,Prodinger,SW,Munarini}).
In particular, Konvalina \cite{Konvalina} considered the number of  $k$-subsets $\{x_1,x_2,\ldots,x_k\}$ of $\mathbb Z_n$ such that
$x_i-x_j\neq 2$ for all $1\leqslant i,j\leqslant k$, and found that the answer is  $\frac{n}{n-k}{n-k\choose k}$ if $n\geqslant 2k+1$.
Hwang \cite{Hwang2} then generalized Konvalina's result to the case $x_i-x_j\neq m$ and deduced that the desired number is
given by the same formula if $n\geqslant mk+1$.
Recently, Mansour and Sun \cite{MS} gave the following unification of Yamamoto's and Hwang's formulas.
\begin{thm}[Mansour-Sun]\label{thm:msun}
Let $m,n, p,k$ be positive integers such that  $n\geqslant mpk+1$. Then the number of $k$-subsets
$\{x_1,x_2,\ldots,x_k\}$ of $\mathbb Z_n$ such that
\begin{align}
x_i-x_j\notin\{m,2m,\ldots,pm\}\quad \text{\rm{(}$1\leqslant i,j\leqslant k${\rm)}},  \label{eq:msun}
\end{align}
is also given by \eqref{eq:yama}.
\end{thm}

A short proof of Theorem \ref{thm:msun}   was given by Guo
\cite{Guo2} by using Rothe's identity. In order to generalize Mansour-Sun's result, Chen, Wang, and
Zhang \cite{CWZ} defined an $m$-AP-block of length $k$ to be a
 sequence  $(x_1,x_2,\ldots,x_k)$ of distinct elements in $\Z_n$
such that  $x_{i+1}-x_i=m$ for $1\leqslant i\leqslant k-1$ and
studied the problem of partitioning $\Z_n$ into $m$-AP-blocks.
The
type of such a partition is defined to be  the type of the multiset of
the lengths of the blocks. For example, the following is a $3$-AP-partition of
$\Z_{20}$ of type $1^4 2^3 3^2 4^1$:
\begin{align*}
(2),\ (4,7),\ (5,8),\ (6),\ (9,12,15),\ (10),\ (11),\ (13,16,19),\ (14,17,0,3),\ (18,1).
\end{align*}
We need to emphasize that
$(x,x+m,\ldots,x+(n-1)m)$ and $(x+m,x+2m,\ldots,x+(n-1)m,x)$ and so on are deemed as different $m$-AP-blocks
in $\Z_{mn}$. For example, all the 2-AP-partitions of $\Z_6$ of type $3^2$ are
$$
\{(i,i+2,i+4),\:(j+1,j+3,j+5)\}_{i,j=0,2,4}.
$$
Chen, Wang, and Zhang \cite{CWZ} constructed a bijection between $m$-AP-partitions and
$m'$-AP-partitions of $\Z_n$ under some technical conditions, and established the following theorem.
\begin{thm}[Chen-Wang-Zhang]\label{thm:cwz}
Let $m,n,k_1,k_2,\ldots,k_r$ and $i_2,\ldots,i_r$ be positive integers such that $1<i_2<\cdots<i_r$ and
\begin{align}
k_1>(k_2+\cdots+k_r)\big((m-1)(i_r-1)-1\big). \label{eq:cwz}
\end{align}
Then the number of partitions of $\Z_n$ into $m$-AP-blocks of type $1^{k_1}i_2^{k_2}\cdots i_r^{k_r}$
does not depend on $m$, and is given by the cyclic multinomial coefficient
\begin{align}
\frac{n}{k_1+\cdots+k_r}{k_1+\cdots+k_r\choose k_1,\ldots,k_r}. \label{eq:multi}
\end{align}
\end{thm}
If we specialize the type to $1^{n-(p+1)k}(p+1)^k$, then the condition \eqref{eq:cwz} becomes $n\geqslant mpk+1$.
Furthermore, if
$(x_1,x_1+m,\ldots,x_1+pm),\ldots, (x_k,x_k+m,\ldots,x_k+pm)$ are
 the $k$ blocks of length $p+1$ in an $m$-AP-partition of $\Z_n$ of type
$1^{n-(p+1)k}(p+1)^k$, then
the set $\{x_1,\ldots,x_k\}$ satisfies \eqref{eq:msun}, and vice versa.
Therefore Theorem \ref{thm:cwz} implies  Theorem \ref{thm:msun}.

In this paper we shall improve and complete Theorem \ref{thm:cwz} by establishing the following two theorems.
\begin{thm}\label{thm:generalcwz}
Let $m,n$ be positive integers, and let $k_1,k_2,\ldots,k_r$ be nonnegative
 integers such that
$n=k_1+2k_2+\cdots+rk_r$. Let $d=\gcd(m,n)$. If
\begin{align}
\Delta:=n-d(n-k_1-\cdots-k_r)>0,  \label{eq:delta}
\end{align}
then the number of partitions of $\Z_n$ into $m$-AP-blocks of type $1^{k_1}2^{k_2}\cdots r^{k_r}$
is given by \eqref{eq:multi}.
\end{thm}
It is not hard to see that the condition \eqref{eq:delta} is weaker than \eqref{eq:cwz}, i.e., the condition \eqref{eq:cwz}
implies that \eqref{eq:delta}. In other words, for fixed
$n$ and a given type, there are in general many more $m$'s satisfying \eqref{eq:delta} than satisfying \eqref{eq:cwz}. For example,
by Theorem \ref{thm:generalcwz}, the numbers of $m$-AP-partitions of $\Z_{120}$ of type $1^{89}2^3 3^2 5^1 7^2$
are all equal for
\begin{align*}
m=&1,2,3,4,5,7,9,11,13,14,17,19,21,22,23,25,26,27,28,29,31, \\
&33,34,35,37,38,39,41,43,44,46,47,49,51,52,53,55,57,58,59,
\end{align*}
i.e., for $d=1,2,3,4,5$.
However, Theorem \ref{thm:cwz} only asserts that these numbers for $m=1,2,3$ are equal.

\begin{thm}\label{thm:gen2}
Let $k_1,k_2,\ldots,k_r,m,n,d$ and $\Delta$ be given as in Theorem \ref{thm:generalcwz}.
Then the number of partitions of $\Z_n$ into $m$-AP-blocks of type $1^{k_1}2^{k_2}\cdots r^{k_r}$ is
given by
\begin{align*}
\begin{cases}
\displaystyle\frac{n}{k_1+\cdots+k_r}{k_1+\cdots+k_r\choose k_1,\ldots,k_r}
+\frac{n(-1)^{k_2+\cdots+k_r}}{k_2+\cdots+k_r}{k_2+\cdots+k_r\choose k_2,\ldots,k_r},
&\text{if $\Delta=0$,}   \\[15pt]
\displaystyle\frac{n}{k_1+\cdots+k_r}{k_1+\cdots+k_r\choose k_1,\ldots,k_r}
+\clubsuit(-1)^{k_2+\cdots+k_r}{k_2+\cdots+k_r\choose k_2,\ldots,k_r}, &\text{if $\Delta=-d$, }
\end{cases}
\end{align*}
where
$$
\clubsuit=
\begin{cases}n,&\text{if $k_2=0$,}\\[5pt]
\displaystyle n\left(1-\frac{n(1-d^{-1})k_2}{(k_2+\cdots+k_r)(k_2+\cdots+k_r-1)}\right),&\text{if $k_2>0$.}
\end{cases}
$$
\end{thm}

When   the type in Theorem \ref{thm:gen2} is
$1^{n-(p+1)k}(p+1)^k$ again, then $\Delta=n-mpk$.
To assure that
there is an $m$-AP-block of length $p+1$ in $\Z_n$, we need to assume that $n>pm$, which is equivalent to
$k\geqslant 2$ if $\Delta=0$ and $pk>p+1$ if $\Delta=-m$.
As mentioned after Theorem \ref{thm:cwz},
 each family of $k$ $m$-AP-blocks in $\Z_n$
is in one-to-one correspondence with a $k$-subset of $\Z_n$ satisfying
\eqref{eq:msun}, we derive the following two results, which can be viewed as
complements to Theorem \ref{thm:msun}.
\begin{cor}\label{cor:1}
Let $m,p\geqslant 1$, $k\geqslant 2$ and $n=mpk$. Then the number of $k$-subsets
$\{x_1,x_2,\ldots,x_k\}$ of $\mathbb Z_n$ such that
$x_i-x_j\notin\{m,2m,\ldots,pm\}$ for all $1\leqslant i,j\leqslant k$, is
given by
\begin{align}
\frac{n}{n-pk}{n-pk\choose k}+(-1)^{k}\frac{n}{k}.  \label{eq:cor1}
\end{align}
\end{cor}

Actually  the above formula is deduced for $m\geqslant 2$, i.e.,
$n\geqslant (p+1)k$, but it also holds for $m=1$ if we take the convention
$$
\lim_{x\to 0}\frac{n}{x}{x\choose k}=(-1)^{k-1}\frac{n}{k},
$$
and so \eqref{eq:cor1} is equal to $0$ in this case.
Here is  an example for Corollary \ref{cor:1}. For $m=p=k=2$, the number of $2$-subsets $\{x_1,x_2\}$ of $\Z_8$
such that $x_1-x_2,x_2-x_1\notin\{2,4\}$ is equal to
$$
\frac{8}{4}{4\choose 2}+4=16,
$$
and the corresponding subsets are $\{i,i+1\}$ and $\{i,i+3\}$, where $i\in\Z_8$.

\begin{cor}\label{cor:2}
Let $m,p,k\geqslant 1$ with $pk>p+1$ and let $n=mpk-m$. Then the number of $k$-subsets
$\{x_1,x_2,\ldots,x_k\}$ of $\mathbb Z_n$ such that
$x_i-x_j\notin\{m,2m,\ldots,pm\}$ for all $1\leqslant i,j\leqslant k$, is
given by
$$
\begin{cases}
\displaystyle\frac{n}{n-k}{n-k\choose k}+(-1)^{k-1}n(m-2), &\text{if $p=1$},\\[15pt]
\displaystyle\frac{n}{n-pk}{n-pk\choose k}+(-1)^{k}n, &\text{if $p\geqslant 2$.}
\end{cases}
$$
\end{cor}

Similarly, although the above formula is deduced for $mpk-m\geqslant (p+1)k$,
it also holds without this condition. The details are left to
the interested reader.

\medskip
\noindent{\it  Remark.}
For $0<m<n$, let $g_m(n,k)$ denote the number of $k$-subsets
$\{x_1,x_2,\ldots,x_k\}$ of $\mathbb Z_n$ such that
$x_i-x_j\neq m$ for all $1\leqslant i,j\leqslant k$. Hwang \cite[Corollary 2]{Hwang2} obtained
\begin{align}
g_m(n,k)=\sum_{j=0}^{\lfloor d/2\rfloor}(-1)^{nj/d}{d\choose j}\frac{n-2nj/d}{n-k-nj/d}{n-k-nj/d\choose k-nj/d},
\label{eq:hwang}
\end{align}
where $d=\gcd(m,n)$.
Letting $n=mk$ or $n=mk-m$ in \eqref{eq:hwang},
we are led to the $p=1$ case of Corollaries \ref{cor:1} or \ref{cor:2}. However, since there are two cases in
Corollary \ref{cor:2}, it seems impossible to give a formula like \eqref{eq:hwang} to unify
Corollaries \ref{cor:1} and \ref{cor:2} for general $p$.

 We  recall and establish some necessary lemmas in Section 2 and
 prove Theorems \ref{thm:generalcwz} and \ref{thm:gen2} in
Sections 3 and 4, respectively.
Our main idea is the following:
Lemma 2.4 permits us to  reduce the general $m$-AP-partition problem of $\Z_n$ to the case where
$m$ divides $n$. For the latter we  may
write the partition number as a multiple sum, which can be   computed by
applying Raney-Mohanty's identity.

\section{Some lemmas}
\begin{figure}[h]
\begin{center}
\begin{pspicture}(5,3)(7.6,0)
\psset{unit=6pt}%
\pscircle[linewidth=1pt](30,0){10}%
\pscircle*(30,10){.4}%
\pscircle*(33.09,9.51){.4}%
\pscircle*(35.88,8.09){.4}%
\pscircle*(38.09,5.88){.4}%
\pscircle*(39.51,3.09){.4}%
\pscircle*(40,0){.4}%
\pscircle*(39.51,-3.09){.4}%
\pscircle*(38.09,-5.88){.4}%
\pscircle*(35.88,-8.09){.4}%
\pscircle*(33.09,-9.51){.4}%
\pscircle*(30,-10){.4}%
\pscircle*(26.91,9.51){.4}%
\pscircle*(24.12,8.09){.4}%
\pscircle*(21.91,5.88){.4}%
\pscircle*(20.49,3.09){.4}%
\pscircle*(20,0){.4}%
\pscircle*(20.49,-3.09){.4}%
\pscircle*(21.91,-5.88){.4}%
\pscircle*(24.12,-8.09){.4}%
\pscircle*(26.91,-9.51){.4}%
\psline[linewidth=.5pt](36.01,6.01)(38.132,8.132)%
\psline[linewidth=.5pt](33.859,7.574)(35.221,10.247)%
\psline[linewidth=.5pt](38.395,-1.33)(41.358,-1.799)%
\psline[linewidth=.5pt](36.01,-6.01)(38.132,-8.132)%
\psline[linewidth=.5pt](31.33,-8.395)(31.799,-11.358)%
\psline[linewidth=.5pt](22.426,3.859)(19.753,5.221)%
\psline[linewidth=.5pt](23.99,6.01)(21.868,8.132)%
\psline[linewidth=.5pt](21.605,-1.33)(18.642,-1.799)%
\psline[linewidth=.5pt](22.426,-3.859)(19.753,-5.221)%
\psline[linewidth=.5pt](28.67,-8.395)(28.201,-11.358)%
\put(29.5,11.5){\small$0$}%
\put(33,11){\small$1$}%
\put(36.2,9.5){\small$2$}%
\put(39,6.5){\small$3$}%
\put(40.6,3){\small$4$}%
\put(41.1,-.5){\small$5$}%
\put(40.5,-4.2){\small$6$}%
\put(39,-7.7){\small$7$}%
\put(36.3,-10.4){\small$8$}%
\put(33,-12.1){\small$9$}%
\put(29,-12.8){\small$10$}%
\put(25.5,-12.3){\small$11$}%
\put(22,-10.5){\small$12$}%
\put(19.2,-7.8){\small$13$}%
\put(17.5,-4){\small$14$}%
\put(17,-.5){\small$15$}%
\put(17.5,3){\small$16$}%
\put(19.2,6.5){\small$17$}%
\put(22,9.5){\small$18$}%
\put(25.5,11){\small$19$}%
\end{pspicture}
\end{center}\vspace{60pt}
\caption{A $20$-cycle dissection of type $1^4 2^3 3^2 4^1$.\label{Fig1}}
\end{figure}
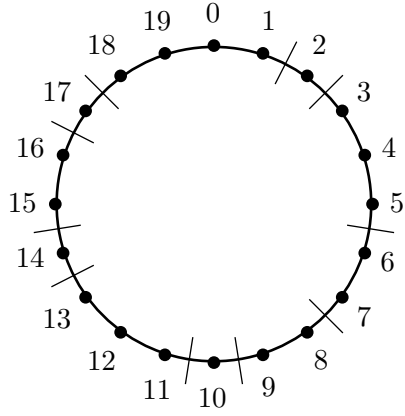

A {\it dissection of an
$n$-cycle} is a $1$-AP-partition  of  $\Z_n$, which
 can be depicted by inserting a bar between any two consecutive blocks on an $n$-cycle.
For example, Figure~\ref{Fig1} illustrates a $20$-cycle
dissection of type $1^4 2^3 3^2 4^1$.
It is easy to see that the number of dissections of $\Z_n$ is given by \eqref{eq:multi}.
Indeed, deleting the segment containing $0$  in any dissection of $n$-cycle of
 type $1^{k_1}2^{k_2}\cdots r^{k_r}$
yields
a dissection of a $(n-i)$-line of type $1^{k_1}\ldots i^{k_i-1}\ldots r^{k_r}$ if
the segment containing $0$ is of length $i$ ($1\leqslant i\leqslant n$). So
the number of such dissections of $n$-cycle is equal to
\begin{align}
i{k_1+\cdots +(k_i-1)+\cdots +k_r\choose k_1,\ldots,k_i-1,\ldots, k_r}. \label{eq:line}
\end{align}
Summing \eqref{eq:line} over all $i$ yields the following known result (see \cite[Lemma 3.1]{CLY}).
\begin{lem}[Chen-Lih-Yeh]\label{lem:cly}
For an $n$-cycle, the number of dissections of type $1^{k_1}2^{k_2}\cdots r^{k_r}$
is given by the cyclic multinomial coefficient \eqref{eq:multi}.
\end{lem}

For any variable $x$ and   nonnegative integers  $k_1,\ldots,k_r$  define  the multinomial coefficient
$$
{x\choose k_1,k_2,\ldots,k_r}:=\frac{x(x-1)\cdots(x-k_1-\cdots-k_r+1)}{k_1!k_2!\cdots k_r!}.
$$
Note that
when $x=k_1+k_2+\ldots+k_r$ the above definition coincides with the classical definition
of multinomial coefficient and
$$
{k_1+\cdots+k_r\choose k_1,\ldots,k_r}={k_1+\cdots+k_r\choose k_2,\ldots,k_r}.
$$

The following convolution formula for multinomial coefficients is due to Raney-Mohanty \cite{Raney,Mohanty66}.
For other proofs of \eqref{eq:rm}, we refer the reader to \cite{Guo3,Strehl,Zeng}.
\begin{lem}[Raney-Mohanty's identity] For any variables $x,y,z_1,\ldots,z_m$ and nonnegative integers $N_1,\ldots,N_m$, there holds
\begin{align}
&\sum_{\substack{0\leqslant t_i\leqslant N_i\\i=1,\ldots,m}}\frac{x}{x-t_1z_1-\cdots-t_m z_m}
{x-t_1z_1-\cdots-t_m z_m\choose t_1,\ldots,t_m} \nonumber\\
&\quad\times \frac{y}{y-(N_1-t_1)z_1-\cdots -(N_m-t_m)z_m}{y-(N_1-t_1)z_1-\cdots -(N_m-t_m)z_m\choose N_1-t_1,\ldots,N_m-t_m} \nonumber\\
&=\frac{x+y}{x+y-N_1z_1-\cdots-N_m z_m}{x+y-N_1z_1-\cdots-N_m z_m\choose N_1,\ldots,N_m}. \label{eq:rm}
\end{align}
\end{lem}

We also need the following elementary   arithmetical result (see \cite[Theorem 5.32 and  Exercise 16 on page 127]{Apostol} or
\cite{Guo2}).
\begin{lem}\label{lem:two}
Let $m, n$ be positive integers.
If $\gcd(m,n)=d$, then there exists an  integer $a$ such that $\gcd(a,n)=1$ and $am\equiv d\pmod{n}$.
\end{lem}

The following is our key lemma.
\begin{lem}\label{lem:three}
If $m,n\geqslant 1$ and $\gcd(m,n)=d$, then  there is a bijection from the set of $m$-AP-partitions of $\Z_n$ to
the set of $d$-AP-partitions of $\Z_n$. Moreover this bijection keeps  the  type of partitions.
\end{lem}
\pf
By Lemma \ref{lem:two}, there exists an inversible element  $a\in \Z_n$ such that
$am=d$.  Let  $a^{-1}$  be the  inverse of $a$.
For any subset $B$ of $\Z_n$ and $x\in \Z_n$,  let $xB=\{xb\colon b\in B\}$.
If $\{B_1, B_2, \ldots, B_s\}$ is  an $m$-AP-partition of $\Z_n$,
then $\{aB_1, aB_2, \ldots, aB_s\}$ is a $d$-AP-partition of $\Z_n$.
Conversely, if
 $\{C_1, C_2, \ldots, C_s\}$  is a $d$-AP-partition of $\Z_n$, then
 $\{a^{-1}C_1, a^{-1}C_2, \ldots,  a^{-1}C_s\}$ is an $m$-AP-partition of $\Z_n$.
Obviously, this correspondence keeps  the  type of partitions.
This proves the lemma.
\qed

It follows from Lemma \ref{lem:three} that if there exists an $m$-AP-partition of $\Z_n$ of
a given type $1^{k_1}2^{k_2}\cdots r^{k_r}$ ($k_r>0$) then $\gcd(m,n)r\leqslant n$.
\section{Proof of Theorem~\ref{thm:generalcwz} }
By Lemma \ref{lem:three}, it suffices to consider the case where $m$ divides $n$, i.e., $d=m$.
Let $n=mn_1$ and divide $\Z_n$ into $m$ subsets of the same cardinality $n_1$:
$$
\mathbb Z_{n,j}=\{mi+j\colon i=0,\ldots,n_1-1\},\qquad 0\leqslant j\leqslant m-1.
$$
Hence $\mathbb Z_n=\biguplus_{j=0}^{m-1}\mathbb Z_{n,j}$.
Let ${ \cal B}=\{B_1, B_2\,\ldots, B_s\}$ be an $m$-AP-partition of $\Z_n$  of
type $1^{k_1}2^{k_2}\cdots r^{k_{r}}$ ($r\leqslant n_1$).
Then $B_{i,j}=\Z_{n,j}\cap B_i$ is equal to $\emptyset$ or  $B_i$
for $1\leqslant i\leqslant s$ and $0\leqslant j\leqslant m-1$.
Furthermore, since the  transformation $x\mapsto (x-j)/m$ maps
each  $m$-AP-block $B_{i,j}$  of $\Z_{n,j}$ ($0\leqslant j\leqslant m-1$)  to a 1-AP-block  $B'_{i,j}$ of $\Z_{n_1}$,
each $m$-AP-partition ${ \cal B}_j=\{B_{1,j}, \ldots, B_{s,j}\}$ corresponds bijectively to a 1-AP-partition  ${\cal B}'_{j}$ of $\Z_{n_1}$
with  the same type.
Thus, we have established a bijection between the set of $m$-AP-partitions  of
 $\Z_n$ and the set of $m$-tuples of
  $1$-AP-partitions of $\Z_{n_1}$:
  ${\cal B}\leftrightarrow ({\cal B}'_{0}, \ldots, {\cal B}'_{m-1})$.

Now assume that  the $m$-AP-partition ${\cal B}$ is of
type $1^{k_1}2^{k_2}\cdots r^{k_{r}}$ ($r\leqslant n_1$), and the corresponding
$1$-AP-partition ${\cal B}'_{j}$ is of type $1^{k_{1,j}}2^{k_{2,j}}\cdots r^{k_{r,j}}$  ($0\leqslant j\leqslant m-1$).
Clearly,
\begin{align}
\begin{cases}
k_{2,0}+k_{2,1}+\cdots+k_{2,m-1}=k_2,\\
k_{3,0}+k_{3,1}+\cdots+k_{3,m-1}=k_3,\\
\cdots \\
k_{r,0}+k_{r,1}+\cdots+k_{r,m-1}=k_r .
\end{cases}  \label{eq:eqs}
\end{align}
By Lemma \ref{lem:cly} and noticing that
$n_1=k_{1,j}+2k_{2,j}+\cdots+rk_{r,j}$, the number of $1$-AP-partitions of $\Z_{n_1}$  of type
$1^{k_{1,j}}2^{k_{2,j}}\cdots r^{k_{r,j}}$
is equal to
\begin{align*}
&\hskip -3mm
\frac{n_1}{k_{1,j}+k_{2,j}+\cdots+k_{r,j}}{k_{1,j}+k_{2,j}+\cdots+k_{r,j}\choose k_{1,j},k_{2,j},\ldots,k_{r,j}}\\
&=
\frac{n_1}{n_1-k_{2,j}-\cdots-(r-1)k_{r,j}}{n_1-k_{2,j}-\cdots-(r-1)k_{r,j}
\choose k_{2,j},\ldots,k_{r,j}}.
\end{align*}
For $m,n\geqslant 1$ let
$f_{m,n}(k_1,\ldots, k_r)$ be the number of partitions of $\Z_n$ into $m$-AP-blocks of type $1^{k_1}2^{k_2}\cdots r^{k_r}$.
Then
\begin{align}
f_{m,n}(k_1,\ldots, k_r)=\sum_{(k_{i,j})}\prod_{j=0}^{m-1}
\frac{n_1}{n_1-k_{2,j}-\cdots-(r-1)k_{r,j}}{n_1-k_{2,j}-\cdots-(r-1)k_{r,j}\choose
k_{2,j},\ldots,k_{r,j}}, \label{eq:fmn}
\end{align}
where the summation is over all matrices $(k_{i,j})_{\substack{2\leqslant i\leqslant r\\ 0\leqslant j\leqslant m-1}}$
of nonnegative integral coefficients
$k_{i,j}$ satisfying
\eqref{eq:eqs} and
\begin{align}
\begin{cases}
n_1-k_{2,0}-\cdots-(r-1)k_{r,0}>0,\\
n_1-k_{2,1}-\cdots-(r-1)k_{r,1}>0,\\
\cdots \\
n_1-k_{2,m-1}-\cdots-(r-1)k_{r,m-1}>0 .
\end{cases}  \label{eq:condnew}
\end{align}

Recall that
\begin{align*}
\Delta &=n-m(n-k_1-\cdots-k_r)
=mn_1-m(k_2+\cdots+(r-1)k_r).
\end{align*}
If $\Delta>0$, then we have
$
n_1>k_2+\cdots+(r-1)k_r,
$
and thus all nonnegative integral solutions to \eqref{eq:eqs} also satisfy \eqref{eq:condnew}
as $k_i\geqslant k_{i,j}$ ($2\leqslant i\leqslant r$, $0\leqslant j\leqslant m-1$).

It remains to prove that the right-hand side of \eqref{eq:fmn} is  equal to \eqref{eq:multi},
namely
\begin{align}\label{eq5}
\frac{mn_1}{mn_1-k_{2}-\cdots-(r-1)k_{r}}{mn_1-k_{2}-\cdots-(r-1)k_{r}\choose k_{2},\ldots,k_{r}}.
\end{align}
We proceed by induction on $m\geqslant 1$.  This is equivalent to repeatedly applying
Raney-Mohanty's identity \eqref{eq:rm}.
The case $m=1$ is obviously true. Suppose that the formula is true for $m-1$ with $m\geqslant 2$ and let
$k_{i,0}+k_{i,1}+\cdots + k_{i,m-2}=k_i'$ be fixed for $i=2,\ldots, r$. Then
\begin{align*}
&\sum_{k_{i,0}, \ldots, k_{i,m-2}\atop i=2,\ldots, r}\prod_{j=0}^{m-2}
\frac{n_1}{n_1-k_{2,j}-\cdots-(r-1)k_{r,j}}{n_1-k_{2,j}-\cdots-(r-1)k_{r,j}\choose
k_{2,j},\ldots,k_{r,j}}\\
&=\frac{(m-1)n_1}{(m-1)n_1-k'_{2}-\cdots-(r-1)k'_{r}}{(m-1)n_1-k'_{2}-\cdots-(r-1)k'_{r}\choose k'_{2},\ldots,k'_{r}}.
\end{align*}
Plugging this into \eqref{eq:fmn}  yields
\begin{align*}
&\hskip -3mm
f_{m,n}(k_1,\ldots, k_r)\\
&=\sum_{k_i'+k_{i,m-1}=k_i\atop i=2,\ldots, r}\frac{(m-1)n_1}{(m-1)n_1-k'_{2}-\cdots-(r-1)k'_{r}}{(m-1)n_1-k'_{2}-\cdots-(r-1)k'_{r}\choose k'_{2},\ldots,k'_{r}}\\
&\times \frac{n_1}{n_1-k_{2,m-1}-\cdots-(r-1)k_{r,m-1}}{n_1-k_{2,m-1}-\cdots-(r-1)k_{r,m-1}\choose
k_{2,m-1},\ldots,k_{r,m-1}}
\end{align*}
which is \eqref{eq5} by applying  Raney-Mohanty's identity \eqref{eq:rm}.

\section{Proof of Theorem~\ref{thm:gen2}}
For the case $\Delta=0$ or $\Delta=-m$, the number $f_{m,n}(k_1,\ldots, k_r)$ is again given by \eqref{eq:fmn}.
However, we will meet with
$
n_1-k_{2,j}-\cdots-(r-1)k_{r,j}\leqslant 0
$
for some $0\leqslant j\leqslant m-1$ in some nonnegative integral solutions $(k_{i,j})$ to \eqref{eq:eqs}.
It is convenient here to consider a more general form of \eqref{eq:fmn} as follows. For any variable
$x$, let
$f_{m,n}(x;k_1,\ldots, k_r)$ be the following expression
\begin{align*}
\sum_{(k_{i,j})}\prod_{j=0}^{m-1}
\frac{x}{x-k_{2,j}-\cdots-(r-1)k_{r,j}}{x-k_{2,j}-\cdots-(r-1)k_{r,j}\choose k_{2,j},\ldots,k_{r,j}},
\end{align*}
where $(k_{i,j})$ ranges over the same integral matrices as \eqref{eq:fmn}.

Let $M$ be the set of all nonnegative integral matrices $(k_{i,j})_{\substack{2\leqslant i\leqslant r\\ 0\leqslant j\leqslant m-1}}$
satisfying \eqref{eq:eqs}, and let
$S$ be the set of all $(k_{i,j})$ in $M$ such that \eqref{eq:condnew} does not hold. Then
\begin{align}
&\hskip -3mm
f_{m,n}(x;k_1,\ldots, k_r) \nonumber\\
&=\sum_{(k_{i,j})\in M}\prod_{j=0}^{m-1}
\frac{x}{x-k_{2,j}-\cdots-(r-1)k_{r,j}}{x-k_{2,j}-\cdots-(r-1)k_{r,j}\choose k_{2,j},\ldots,k_{r,j}} \nonumber\\
&\quad{}-\sum_{(k_{i,j})\in S}\prod_{j=0}^{m-1}
\frac{x}{x-k_{2,j}-\cdots-(r-1)k_{r,j}}{x-k_{2,j}-\cdots-(r-1)k_{r,j}\choose k_{2,j},\ldots,k_{r,j}}.  \label{eq:fffmn}
\end{align}

When $\Delta=0$, we have $ n_1=k_2+2k_3+\cdots+(r-1)k_r$, and $S$ reduces to
$$
S_1:=\left\{(k_{i,j})\colon\
 \text{for some $j_0$ and all $i$, we have $k_{i,j_0}=k_{i}$ and $k_{i,j}=0$ if $j\neq j_0$}\right\}.
$$
So the second summation on the right-hand side of \eqref{eq:fffmn} becomes
$$
\frac{mx}{x-n_1}{x-n_1\choose k_{2},\ldots,k_{r}},
$$
while
 the first summation can be summed by using  Raney-Mohanty's identity.
It follows that
\begin{align}
&\hskip -3mm
f_{m,n}(x;k_1,\ldots, k_r) \nonumber \\
&=\frac{mx}{mx-k_{2}-\cdots-(r-1)k_{r}}{mx-k_{2}-\cdots-(r-1)k_{r}\choose k_{2},\ldots,k_{r}}
-\frac{mx}{x-n_1}{x-n_1\choose k_{2},\ldots,k_{r}}.  \label{eq:fin1}
\end{align}
Letting $x=n_1$ in \eqref{eq:fin1} and noticing the following fact
\begin{align}
\lim_{z\to 0}\frac{1}{z}{z\choose a_1,\ldots,a_{s}}
=\frac{(-1)^{a_1+\cdots+a_{s}-1}}{a_1+\cdots+a_{s}}{a_1+\cdots+a_{s}\choose a_1,\ldots,a_{s}}, \label{eq:limit0}
\end{align}
one obtains the first formula in Theorem~\ref{thm:gen2}.

When $\Delta=-m$, we have
$
n_1=k_2+2k_3+\cdots+(r-1)k_r-1.
$
If $k_2=0$, then $S=S_1$, while if $k_2>0$, then
\begin{align*}
S&=S_1
\cup
 \{(k_{i,j})\colon\
 \text{for some $j_0\neq j_1$, we have $k_{2,j_0}=k_{2}-1$, $k_{2,j_1}=1$}, \\
 &\hskip 4cm \text{$k_{i,j_0}=k_{i}$ ($2<i\leqslant r$)
 and $k_{i,j}=0$ otherwise}\}.
\end{align*}
It follows that
\begin{align}
&\hskip -3mm
f_{m,n}(x;k_1,\ldots, k_r) \nonumber \\
&=\frac{mx}{mx-k_{2}-\cdots-(r-1)k_{r}}{mx-k_{2}-\cdots-(r-1)k_{r}\choose k_{2},\ldots,k_{r}} \nonumber \\
&\quad{}-\frac{mx}{x-n_1-1}{x-n_1-1\choose k_{2},\ldots,k_{r}}
-\chi(k_2>0)\frac{m(m-1)x^2}{x-n_1}{x-n_1\choose k_{2}-1,k_3,\ldots,k_{r}}.     \label{eq:k2>0}
\end{align}
Letting $x=n_1$ in \eqref{eq:k2>0} and using \eqref{eq:limit0} and
$$
{-1\choose a_1,\ldots,a_s}=(-1)^{a_1+\cdots+a_s}{a_1+\cdots+a_s\choose a_1,\ldots,a_s},
$$ we obtain the second formula in Theorem~\ref{thm:gen2}.

\medskip
\noindent{\it Remark.} It is also possible to compute
$f_{m,n}(k_1,\ldots, k_r)$ for the case $\Delta=-2\gcd(m,n)$ or
$\Delta=-3\gcd(m,n)$. But the result is more complicated and is
omitted here.

\medskip
\noindent{\bf Acknowledgments.} We thank the referees for helpful comments on a previous version of this paper.
This work was done during the first author's visit to Institut Camille Jordan of Univerit\'e
Lyon I, and was supported by Project MIRA 2007 de la  R\'egion Rh\^one-Alpes.
The first author was also supported by Shanghai Leading Academic Discipline Project, Project Number: B407.

\end{document}